\documentclass{amsart}
\usepackage{amsmath}
\usepackage{amssymb}
\usepackage[english]{babel}
\usepackage{graphics}
\usepackage[all]{xy}
\usepackage{amsrefs}
\usepackage{amsmath}
\usepackage{amsthm}
\usepackage{amssymb}
\usepackage{amsfonts}
\usepackage{amsxtra}
\usepackage{tikz}
\usepackage{multicol}
\usepackage{subfigure}
\usepackage{centernot}
\usepackage{soul}
\usepackage{xcolor}
\usepackage{longtable}
\usepackage{layout}
\usepackage{geometry}
\usepackage{amsthm}

\usepackage{amstext} 
\usepackage{array}   
\newcolumntype{L}{>{$}c<{$}} 
\usepackage{longtable}

\usepackage[colorlinks=true,linkcolor=blue,urlcolor=blue, citecolor=blue]{hyperref}
\usepackage{fullpage}
\usepackage{epsfig}
\usepackage{verbatim}
\usepackage{enumitem}
\usepackage{algorithmicx}
\usepackage{float}
\usepackage{algpseudocode}

\usepackage[justification=centering]{caption}

\newtheorem{Theorem}{Theorem}[section]

\newtheorem{Lemma}[Theorem]{Lemma}

\theoremstyle{definition}

\theoremstyle{remark}

\newcommand{\cod}{\operatorname{cod}}

\newcommand{\irr}{\operatorname{Irr}}
\newcommand{\PSL}{\normalfont{\mbox{PSL}}}
\newcommand{\PSU}{\normalfont{\mbox{PSU}}}
\usepackage{graphicx}
\usepackage{adjustbox}

\makeatletter
\newcommand*\bigcdot{\mathpalette\bigcdot@{.5}}
\newcommand*\bigcdot@[2]{\mathbin{\vcenter{\hbox{\scalebox{#2}{$\m@th#1\bullet$}}}}}
\makeatother

\title{On the Characterization of Alternating Groups by Codegrees}
\author[M.~Dolorfino]{Mallory Dolorfino}
\author[L. ~Martin]{Luke Martin}
\author[Z. ~Slonim]{Zachary Slonim}
\author[Y. ~Sun]{Yuxuan Sun}
\author[Y. ~Yang]{Yong Yang}
\address{Mallory~Dolorfino\\ Kalamazoo College \\ Kalamazoo, Michigan, USA \\
  \href{mailto:mallory.dolorfino19@kzoo.edu}
  {{\ttfamily\upshape mallory.dolorfino19@kzoo.edu}}}

\address{Luke~Martin\\ Gonzaga University\\ Spokane, Washington, USA \\
  \href{mailto:lwmartin2019@gmail.com}
  {{\ttfamily\upshape lwmartin2019@gmail.com}}}

\address{Zachary~Slonim\\ University of California, Berkeley \\ Berkeley, California, USA \\
  \href{mailto:zachslonim@berkeley.edu}
  {{\ttfamily\upshape zachslonim@berkeley.edu}}}

\address{Yuxuan~Sun\\ Haverford College \\ Haverford, Pennsylvania, USA \\
  \href{mailto:ysun1@haverford.edu}
  {{\ttfamily\upshape ysun1@haverford.edu}}}

\address{Yong~Yang\\ Texas State University \\ San Marcos, Texas, USA \\
  \href{mailto:yang@txstate.edu}
  {{\ttfamily\upshape yang@txstate.edu}}}


\begin{document}

\maketitle

\begin{abstract}

    Let $G$ be a finite group and $\mathrm{Irr}(G)$ the set of all irreducible complex characters of $G$. Define the codegree of $\chi \in \mathrm{Irr}(G)$ as $\mathrm{cod}(\chi):=\frac{|G:\mathrm{ker}(\chi) |}{\chi(1)}$ and denote by $\mathrm{cod}(G):=\{\mathrm{cod}(\chi) \mid \chi\in \mathrm{Irr}(G)\}$ the codegree set of $G$. Let $\mathrm{A}_n$ be an alternating group of degree $n \ge 5$. In this paper, we show that $\mathrm{A}_n$ is determined up to isomorphism by $\cod(\mathrm{A}_n)$.

\end{abstract}

\begin{section}{Introduction}

Let $G$ be a finite group and $\mathrm{Irr}(G)$ the set of all irreducible complex characters of $G$. For any $\chi \in \mathrm{Irr}(G),$ define the codegree of $\chi$ as $\mathrm{cod}(\chi) := \frac{|G:\mathrm{ker}(\chi)|}{\chi(1)}.$ Then define the codegree set of $G$ as $\mathrm{cod}(G) := \{\mathrm{cod}(\chi) \mid \chi \in \mathrm{Irr}(G)\}.$ The concept of codegrees was originally considered in \cite{Chillag}, where the codegree was defined as $\mathrm{cod}(\chi) := \frac{|G|}{\chi(1)},$ and it was later modified to its current definition by \cite{Qian} so that $\mathrm{cod}(\chi)$ is the same for $G$ and $G/N$ when $N \leq \mathrm{ker}(\chi).$ Several properties of codegrees have been studied, such as the relationship between the codegrees and element orders, codegrees of $p$-groups, and groups with few codegrees.

The codegree set of a group is closely related to the character degree set of a group, which is defined as $\mathrm{cd}(G) := \{\chi(1) \mid  \chi \in \mathrm{Irr}(G)\}.$ The relationship between the character degree set and a group's structure is an active area of research -- many properties of a group's structure are largely determined by its character degree set. In 1990, Bertram Huppert made the following conjecture about the relationship between a simple group $H$ and a finite group $G$ having equal character degree sets.

{\bf Huppert's Conjecture:} Let $H$ be a finite nonabelian simple group and $G$ a finite group such that $\mathrm{cd}(H) = \mathrm{cd}(G).$ Then $G \cong H \times A,$ where $A$ is an abelian group.

Huppert's conjecture has been verified for many cases such as the alternating groups, sporadic groups, and simple groups of Lie type with low rank, but it has yet to be verified for simple groups of Lie type with high rank. Recently, a similar conjecture related to codegrees has been posed.

{\bf Codegree Version of Huppert's Conjecture:} Let $H$ be a finite nonabelian simple group and $G$ a finite group such that $\mathrm{cod}(H) = \mathrm{cod}(G).$ Then $G \cong H.$

This conjecture appears in the \emph{Kourovka Notebook of Unsolved Problems in Group Theory} as question 20.79 \cite{Khukhro}. It has been verified for $\mathrm{PSL}(2,q)$, $\mathrm{PSL}(3,4),$ $\mathrm{Alt}_7,$ $\mathrm{J}_1$, $^2B_2(2^{2f+1})$ where $f \geq 1$, $\mathrm{M}_{11}, \mathrm{M}_{12}, \mathrm{M}_{22}, \mathrm{M}_{23}$ and $\mathrm{PSL}(3,3)$ by \cites{Ahanjideh,Bahri,Gintz}. The conjecture has also been verified for $\PSL(3,q)$ and $\PSU(3,q)$ in \cites{lmy} and ${}^2G_2(q)$ in \cites{gzy}. Recently, the authors verified the conjecture for all sporadic simple groups in ~\cite{Dolorfino}.

In this paper, we provide a general proof verifying this conjecture for all alternating groups of degree greater than or equal to $5$. The methods used may be generalized to simple groups of Lie type, giving promising results for characterizing all simple groups by their codegree sets.

\begin{Theorem}\label{thm11}

    Let $\mathrm{A}_n$ be an alternating group of degree $n \ge 5$ and $G$ a finite group. If $\mathrm{cod}(G)= \mathrm{cod}(\mathrm{A}_n)$, then $G \cong \mathrm{A}_n$.

\end{Theorem}

Throughout the paper, we follow the notation used in Isaacs' book \cite{Isaacs} and the ATLAS of Finite Groups \cite{Atlas}.

\end{section}

\newpage
\begin{section}{Preliminary Results}\label{theory}

We first introduce some lemmas which will be used later.

\begin{Lemma}\label{oneextend}
    \cite{Moreto}*{Lemma 4.2} Let $S$ be a finite nonabelian simple group. Then there exists $1_S\neq \chi \in \mathrm{Irr}(S)$ that extends to $\mathrm{Aut}(S)$.
\end{Lemma}
\begin{Lemma}\label{prodsimple}
    \cite{James}*{Theorem 4.3.34} Let $N$ be a minimal normal subgroup of $G$ such that $N=S_1\times\dots\times S_t$ where $S_i \cong S$ is a nonabelian simple group for each $i=1,\dots,t$. If $\chi \in \mathrm{Irr}(S)$ extends to $\mathrm{Aut}(S)$, then $\chi\times\dots\times\chi\in\mathrm{Irr}(N)$ extends to $G$.
\end{Lemma}
\begin{Lemma}\label{perfect}
    \cite{Gintz}*{Remark 2.6} Let $G$ be a finite group and $S$ a finite nonabelian simple group with $\mathrm{cod}(G)=\mathrm{cod}(S)$. Then $G$ is a perfect group.
\end{Lemma}

\begin{Lemma}\label{order}
    \cite{Hung} Let $G$ be a finite group and $S$ a finite nonabelian simple group such that $\mathrm{cod}(S) \subseteq$ $\mathrm{cod}(G)$. Then $|S|$ divides $|G|$.
\end{Lemma}

\begin{Lemma}\label{codegreeSubset}
    Let $G$ be a finite group with $N \trianglelefteq G.$ Then $\cod(G/N) \subseteq \cod(G)$.
\end{Lemma}
\begin{proof}
From \cite{Isaacs}*{Lemma 2.22}, we can define $\irr(G/N) = \{ \hat{\chi}(gN) = \chi(g) \mid \chi \in \irr(G) \text{ and } N \subseteq \operatorname{ker}(\chi)\}$. Take any $\hat{\chi} \in \mathrm{Irr}(G/N)$. By definition, we know that $\hat{\chi}(1) = \chi(1),$ so the denominators of $\mathrm{cod}(\hat{\chi})$ and $\mathrm{cod}(\chi)$ are equal. In addition, $\mathrm{ker}(\hat{\chi}) \cong \mathrm{ker}(\chi) / N$, so $|\mathrm{ker}(\chi)| = |N| \cdot |\mathrm{ker}(\hat{\chi})|$. Thus $|G/N : \ker(\hat{\chi})| = \frac{|G|/|N|}{|\ker(\chi)|/|N|} = \frac{|G|}{|\ker(\chi)|},$ so $\mathrm{cod}(\hat{\chi}) = \mathrm{cod}(\chi)$ and therefore $\mathrm{cod}(G/N) \subseteq \mathrm{cod}(G).$
\end{proof}
\begin{Lemma}\label{subset}
    Let $G$ be a finite group with normal subgroups $N$ and $M$ such that $N\leq M$. Then, $\mathrm{cod}(G/M)\subseteq \mathrm{cod}(G/N)$.
\end{Lemma}
\begin{proof}
     By the Third Isomorphism Theorem, we know that $G / M \cong (G/N)/(M/N)$ is a quotient of  $G / N$, and by Lemma \ref{codegreeSubset}, $\mathrm{cod}(G / M) \subseteq \mathrm{cod}(G / N)$.
\end{proof}

\end{section}

\begin{Lemma}\label{orderK}
    Let $S$ be a finite nonabelian simple group and $G$ be a nontrivial finite group with $\mathrm{cod}(G)\subseteq \mathrm{cod}(S)$. Then, $|S| < |G|\cdot |\mathrm{Irr}(G)|$.
\end{Lemma}

\begin{proof}
    We know that for each irreducible character $\chi \in \mathrm{Irr}(S), \chi(1)^2 < |S|$. Because $S$ is simple, if $\chi$ is non-trivial, then $\mathrm{ker}(\chi)=1$, so $\mathrm{cod}(\chi)=\frac{|S|}{\chi(1)} > \sqrt{|S|}$. Then, since $\mathrm{cod}(G)\subseteq \mathrm{cod}(S),$ for each irreducible non-trivial character $\psi \in \mathrm{Irr}(G), \mathrm{cod}(\psi) > \sqrt{|S|}$. Thus, $\frac{|G:\mathrm{ker} (\psi)|}{\psi(1)}>\sqrt{|S|}$ which implies that $\frac{|G|}{|\mathrm{ker}(\psi)|\sqrt{|S|}}>\psi(1).$ So, $\psi(1) < \frac{|G|}{\sqrt{|S|}}$. Then $\sum_{\psi \in \operatorname{Irr}(G)} \psi(1)^2 < |\operatorname{Irr}(G)| \frac{|G|^2}{|S|}$, and by character theorems, we'll have $|G| < |\operatorname{Irr}(G)|\frac{|G|^2}{|S|}$. Thus $|S| < |G|\cdot |\mathrm{Irr}(G)|$.
\end{proof}

\begin{section}{Main Results}

We start with some lemmas which limit the simple groups whose codegree set can be contained in the codegree set of an alternating group.

\begin{Lemma}\label{alternating}
    Let $H$ be an alternating group of degree $m\neq n$, where $m,n \ge 5$. Then $\cod(H) \not \subseteq\cod(\mathrm{A}_n)$.
\end{Lemma}
\begin{proof}
Suppose $\cod(\mathrm{A}_m)\subseteq \cod(\mathrm{A}_n)$. Then, from Lemma \ref{order}, $|\mathrm{A}_m|$ divides $|\mathrm{A}_n|,$ so $m < n$. Let $a_x$ denote the minimal non-trivial codegree of $\mathrm{A}_x$. We show that $a_{n-1} < a_n$ so that $\cod(A_m) \not\subseteq \cod(A_n)$ follows immediately.

We know that irreducible representations of the symmetric group $\mathrm{S}_n$ are in one-to-one correspondence with the partitions of $n$. Let $\lambda$ be a partition of $n$ and $V_\lambda$ be the corresponding irreducible representation of $\mathrm{S}_n$. We note that a partition of $n$ can be visualized by a Young diagram and we let $h_\lambda(i,j)$ be the hook length of the $(i,j)^{th}$ square of the Young diagram corresponding to $\lambda$, i.e. the number of cells $(a,b)$ of $\lambda$ such that $a = i$ and $b \geq j$ or $b = j$ and $a \geq i$. By the hook length formula, $\frac{n!}{\mathrm{dim}(V_\lambda)}= \prod h_\lambda(i,j) := H_\lambda$. Let $U_\lambda$ be an irreducible constituent of the restriction of $V_\lambda$ to $\mathrm{A}_n$, $\operatorname{Res}^{\mathrm{S}_n}_{\mathrm{A}_n} V_\lambda$. If $\lambda$ is not self-conjugate ($\lambda \neq \lambda'$), then $\operatorname{Res}^{\mathrm{S}_n}_{\mathrm{A}_n} V_\lambda$ remains irreducible, so $U_\lambda = \operatorname{Res}^{\mathrm{S}_n}_{\mathrm{A}_n} V_\lambda$. In this case, $\frac{n!}{\mathrm{dim}(U_\lambda)}=H_\lambda$. If $\lambda$ is self-conjugate, then the restriction of $V_\lambda$ to $\mathrm{A}_n$ splits into two irreducible representations of the same dimension, so $\mathrm{dim}(U_\lambda)=\frac{1}{2}\mathrm{dim}(V_\lambda)$. In this case, $\frac{n!}{\mathrm{dim}(U_\lambda)}=2H_\lambda$.




Now, $a_n =\min\{ \frac{n!/2}{\mathrm{dim}(U_\lambda)} \mid U_\lambda \in \mathrm{Irr}(\mathrm{A}_n)\} = \frac{1}{2}\min(\{H_\lambda \ |\ \lambda \neq \lambda'\}\cup\{2H_\lambda \ |\ \lambda = \lambda'\})$. We want to show that $a_{n-1} < a_n$. First, assume that $a_n = \frac{1}{2} 2H_\lambda$ for some $\lambda = \lambda'$. Then we can remove a square from $\lambda$ to give a non-self-conjugate partition $\mu$ of $n-1$. Since $H_\mu < H_\lambda < 2H_\lambda$ and $a_{n-1}\leq \frac{1}{2}H_\mu$, we know $a_{n-1} < a_n$.

Now assume that $a_n = \frac{1}{2}H_\lambda$ for some $\lambda \neq \lambda'$. Then if $n \ge 3,$ we can remove a square from $\lambda$ to obtain a non-self-conjugate partition $\mu$ of $n-1$. Since $H_\mu < H_\lambda$ and $a_{n-1}\leq \frac{1}{2}H_\mu$,  $a_{n-1} < a_n.$ Thus, if $m<n,$ then $a_m<a_n,$ contradicting the assumption that $\cod(\mathrm{A}_m)\subseteq \cod(\mathrm{A}_n).$
\end{proof}

\begin{Lemma}\label{sporadic}
    Let $H$ be a sporadic simple group or the Tits group. Then if $n \ge 5,$ $\cod(H) \not \subseteq\cod(\mathrm{A}_n)$.
\end{Lemma}

\begin{proof}
In search of a contradiction, let $H$ be a sporadic simple group or the Tits group such that $\cod(H) \subseteq\cod(\mathrm{A}_n)$. From Lemmas \ref{orderK} and \ref{order}, we deduce a tight restriction on the order of $H$. Namely, $|H|=|\mathrm{A}_n|/k$ where $1\leq k < |\mathrm{Irr}(H)|$ is an integer. Now, for each sporadic (or Tits) group $H$, we can computationally check (using Julia \cite{Julia}) which alternating groups $\mathrm{A}_n$ satisfy both $|H|$ divides $|\mathrm{A}_n|$ and $\frac{|\mathrm{A}_n|}{|H|}<|\mathrm{Irr}(H)|$. We find only one possible exception: $\mathrm{A}_n=\mathrm{A}_{10}$ and $H=J_2$ where $\frac{|\mathrm{A}_{10}|}{|J_2|}=3<21=|\mathrm{Irr}(J_2)|$. In this case, we check that $\cod(J_2)\not \subseteq \cod(\mathrm{A}_{10})$ using the ATLAS \cite{Atlas}.
\end{proof}

\begin{Lemma}\label{classical}
    Let $H$ be a classical simple group of Lie type. Then $\cod(H) \not \subseteq\cod(\mathrm{A}_n)$ for all $n \ge 5$.
\end{Lemma}
\begin{proof}
There are $6$ families of classical simple groups of Lie type. These are $\mathrm{PSL}(m+1,q), \Omega(2m+1,q), \mathrm{PSp}(2m,q), \mathrm{O}^+(2m,q), \mathrm{PSU}(m+1,q),$ and $\mathrm{O}^-(2m,q)$.l We prove the lemma in each case. Let $k(G)$ denote the number of conjugacy classes of $G$, we reproduce \cite{Fulman}*{Table 2} for reference.

\begin{table}[H]
    \centering
    \caption{Class Numbers for Classical Groups}
    \begin{tabular}{ |c|c|c| }
    \hline
    $G$ & $k(G)\leq$ & Comments \\
    \hline \hline
    $\mathrm{SL}(n,q)$ & $2.5q^{n-1}$ & \\ \hline
    $\mathrm{SU}(n,q)$ & $8.26q^{n-1}$ & \\ \hline
    $\mathrm{Sp}(2n,q)$ & $10.8q^n$ & $q$ odd \\ \hline
    $\mathrm{Sp}(2n,q)$ & $15.2q^n$ & $q$ even \\ \hline
    $\mathrm{SO}(2n+1,q)$ & $7.1q^n$ & $q$ odd \\ \hline
    $\Omega(2n+1,q)$ & $7.3q^n$ & $q$ odd \\ \hline
    $\mathrm{SO}^\pm(2n,q)$ & $7.5q^n$ & $q$ odd \\ \hline
    $\Omega^\pm(2n,q)$ & $6.8q^n$ & $q$ odd \\ \hline
    $\mathrm{O}^\pm(2n,q)$ & $9.5q^n$ & $q$ odd \\ \hline
    $\mathrm{SO}^\pm(2n,q)$ & $14q^n$ & $q$ even \\ \hline
    $\mathrm{O}^\pm(2n,q)$ & $15q^n$ & $q$ even \\ \hline
    \end{tabular}
    \label{tbl1}
\end{table}

\begin{enumerate}
\item Let $H=\mathrm{PSL}(m+1,q)$ where $q=p^k$ and $m\geq 1$. From the order formula found in \cite{Carter}, $q^{m(m+1)/2}$ divides $|\mathrm{PSL}(m+1,q)|$. From Legendre's formula, we know that for any prime $p$, $|n!|_p\leq p^{\frac{n}{p-1}}$. If $q=p^k,$ then we have $|n!|_q\leq q^{\frac{n}{k(p-1)}}$ and thus $|\mathrm{A}_n|_q\leq q^{\frac{n}{k(p-1)}}$. By Lemma \ref{order}, $|\mathrm{PSL}(m+1,q)|$ divides $|\mathrm{A}_n|,$ so $q^{m(m+1)/2}$ divides $|\mathrm{A}_n|.$ Thus $\frac{m(m+1)}{2}\leq \frac{n}{k(p-1)},$ giving $n\geq \frac{m(m+1)k(p-1)}{2}$. Therefore, $|\mathrm{A}_n|\geq \left|\mathrm{A}_{\frac{m(m+1)k(p-1)}{2}}\right|$.

Now, we note that $k(\mathrm{PSL}(m+1,q)) \leq k(\mathrm{SL}(m+1,q))$ since $\mathrm{PSL}(m+1,q)$ is a quotient of $\mathrm{SL}(m+1,q).$ Then from Table \ref{tbl1}, we have that $|\mathrm{Irr}(\mathrm{PSL}(m+1,q))|=k(\mathrm{PSL}(m+1,q)) \leq k(\mathrm{SL}(m+1,q))\leq 2.5q^m$. Applying Lemma \ref{orderK} gives $|\mathrm{A}_n|<|\mathrm{PSL}(m+1,q)|\cdot |\mathrm{Irr}(\mathrm{PSL}(m+1,q))|. $ Hence $|\mathrm{A}_{\frac{m(m+1)k(p-1)}{2}}| < |\mathrm{PSL}(m+1,q)|\cdot2.5q^m$. Now we show that if we consider the left and right sides as functions of $m$ with constants $p$ and $k,$ then asymptotically, the value of $|\mathrm{A}_{\frac{m(m+1)k(p-1)}{2}}|$ grows faster than that of $|\mathrm{PSL}(m+1,q)|\cdot2.5q^m.$ We know that the left function behaves asymptotically as $(m^2)!,$ and using the order formula for $\mathrm{PSL}(m+1,q),$ we know that the right function behaves asymptotically as $q^{f(m)},$ where $f(m)$ is a polynomial with degree $2$. Thus the left function grows faster than the right function since $x!>>c^x$ for any constant $c$ when $x$ is large. Similarly, we can prove this result considering the two sides as functions of $p$ and $k.$

Then, we search for the maximum possible value of $m$ which satisfies the inequality given the smallest possible values of $p$ and $k,$ which are $2$ and $1,$ respectively. We find that $m\leq 6$ and, using a similar process for $p$ and $k$, that $p \leq 17$ and $k \leq 63$. Now, we have limited our search to a finite number of groups which we can check in the same way as for the sporadic groups. From this, we find a small list of exceptions, listed in Table ~\ref{tblpsl}:
\begin{table}[H]
    \centering
    \caption{Exceptions satisfying $|\mathrm{PSL}(m+1,q)|$ divides $|\mathrm{A}_n|$ and $|\mathrm{A}_n|<|\mathrm{PSL}(m+1,q)|\cdot2.5q^m$}
    \begin{tabular}{ |c|c|c| }
    \hline
    $m$ & $q$ & $n$ \\
    \hline \hline
1 & 4 & 5 \\ \hline
1 & 4 & 6 \\ \hline
1 & 8 & 7 \\ \hline
1 & 9 & 6 \\ \hline
1 & 9 & 7 \\ \hline
1 & 5 & 5 \\ \hline
1 & 5 & 6 \\ \hline
1 & 7 & 7 \\ \hline
2 & 4 & 8 \\ \hline
2 & 4 & 9 \\ \hline
3 & 2 & 8 \\ \hline
3 & 2 & 9 \\ \hline
    \end{tabular}
    \label{tblpsl}
\end{table}

Now, all of these exceptions can be found in the ATLAS, and it is routine to check that none of these groups satisfy $\mathrm{cod}(\mathrm{PSL}(m+1,q))\subseteq\mathrm{cod}(\mathrm{A}_n)$ unless $\mathrm{PSL}(m+1,q)\cong \mathrm{A}_n$. Thus, if $\mathrm{PSL}(m+1,q)\not\cong \mathrm{A}_n,$ then $\mathrm{cod}(\mathrm{PSL}(m+1,q)\not \subseteq \mathrm{cod}(\mathrm{A}_n)$.

\item Let $H=\Omega(2m+1,q)$ where $q=p^k$ is odd and $m\geq 2$. Note that when $q=2^k$ is even, $\Omega(2m+1,q)\cong \mathrm{PSp}(2m,q),$ which we deal with in the next case. From \cite{Carter}, $q^{m^2}$ divides $|\Omega(2m+1,q)|$. Thus, using Table \ref{tbl1} similarly to above, $|\mathrm{A}_{m^2k(p-1)}| < |\Omega(2m+1,q)|\cdot7.3q^m$. As above, we computationally check that we get a contradiction if $m>2, p>3,$ or $k>1,$ so $m=2,p=3,$ and $k=1$ is the only possibility. We get the list of exceptions listed in Table ~\ref{tblomega} after checking divisibility.

\begin{table}[H]
    \centering
    \caption{Exceptions satisfying $|\Omega(2m+1,q)|$ divides $|\mathrm{A}_n|$ and $|\mathrm{A}_n|<|\Omega(2m+1,q)|\cdot7.3q^m$}
    \begin{tabular}{ |c|c|c| }
    \hline
    $m$ & $q$ & $n$ \\
    \hline \hline
2 & 3 & 9 \\ \hline
    \end{tabular}
    \label{tblomega}
\end{table}

Again, we check the ATLAS and find that $\mathrm{cod}(\Omega(5,3))\not \subseteq \mathrm{cod}(\mathrm{A}_9)$.

\item Let $H=\mathrm{PSp}(2m,q)$ where $q=p^k$ and $m\geq 3$. From \cite{Carter}, $q^{m^2}$ divides $|\mathrm{PSp}(2m,q)|$. Since $\mathrm{PSp}(2m,q)$ is a quotient of $\mathrm{Sp}(2m,q)$, we have $k(\mathrm{PSp}(2m,q)) \leq k(\mathrm{Sp}(2m,q))$. From Table \ref{tbl1}, $|\mathrm{A}_{m^2k(p-1)}| < |\mathrm{PSp}(2m,q)|\cdot15.2q^m$. We computationally check that we get a contradiction if $m>4, p>2,$ or $k>2,$ so $m=3$ or $4,$ $p=2,$ and $k=1$ or $2$ are the only possibilities. We get no exceptions after checking divisibility.

\item Let $H=\mathrm{O}^+(2m,q)$ where $q=p^k$ and $m\geq 4$. From \cite{Carter}, $q^{m(m-1)}$ divides $|O^+(2m,q)|$. Using Table \ref{tbl1}, we have that $|\mathrm{A}_{m(m-1)k(p-1)}| < |O^+(2m,q)|\cdot15q^m$. As above, we computationally check that we get a contradiction if $m>4, p>2,$ or $k>1$ so $m=4, p=2,$ and $k=1$ is the only possibility, and we get no possible exceptions after checking divisibility.

\item Let $H=\mathrm{PSU}(m+1,q)$ where $q=p^k$ and $m\geq 2$. From \cite{Carter}, $q^{m(m+1)/2}$ divides $|\mathrm{PSU}(m+1,q)|$. Since $\mathrm{PSU}(m+1,q)$ is a quotient of $\mathrm{SU}(m+1,q)$, we have $k(\mathrm{PSU}(m+1,q)) \leq k(\mathrm{SU}(m+1,q))$. From Table \ref{tbl1}, $|\mathrm{A}_{\frac{m(m+1)k(p-1)}{2}}| < |\mathrm{PSU}(m+1,q)|\cdot8.26q^m$. Again, we computationally check that we get a contradiction if $m>6, p>7,$ or $k>42$ so $m\leq6, p\leq7,$ and $k\leq 42$ are the only possibilities. We get Table ~\ref{tblpsu} after checking divisibility:

\begin{table}[H]
    \centering
    \caption{Exceptions satisfying $|\mathrm{PSU}(m+1,q)|$ divides $|\mathrm{A}_n|$ and $|\mathrm{A}_n|<|\mathrm{PSU}(m+1,q)|\cdot8.26q^m$}
    \begin{tabular}{ |c|c|c| }
    \hline
    $m$ & $q$ & $n$ \\
    \hline \hline
2 & 3 & 9 \\ \hline
3 & 2 & 9 \\ \hline
    \end{tabular}
    \label{tblpsu}
\end{table}

We check the ATLAS to find that $\mathrm{cod}(\mathrm{PSU}(3,3))\not \subseteq \mathrm{cod}(\mathrm{A}_9),$ and we note that $\mathrm{PSU}(4,2)\cong \Omega(5,3),$ which we have already ruled out.

\item Let $H=\mathrm{O}^-(2m,q)$ where $q=p^k$ and $m\geq 4$. From \cite{Carter}, $q^{m(m-1)}$ divides $|\mathrm{O}^-(2m,q)|$. Thus, using Table \ref{tbl1} similarly to above, $|\mathrm{A}_{m(m-1)k(p-1)}| < |\mathrm{O}^-(2m,q)|\cdot15q^m$. Again, we computationally check that we get a contradiction if $m>5, p>3,$ or $k>3$ so $m\leq5, p\leq3,$ and $k\leq 3$ are the only possibilities, and we get no possible exceptions after checking divisibility.

\end{enumerate}
\end{proof}

\begin{Lemma}\label{exceptional}
    Let $H$ be an exceptional simple group of Lie type. Then if $n \ge 5$, $\cod(H) \not \subseteq\cod(\mathrm{A}_n)$.
\end{Lemma}
\begin{proof}
There are $10$ familes of exceptional simple groups of Lie type (other than the Tits group). These are $E_6(q), E_7(q), E_8(q), F_4(q), G_2(q), ^2E_6(q), ^3D_4(q), ^2B_2(q), ^2F_4(q),$ and $^2G_2(q)$. We prove the lemma in each case. First, we reproduce \cite{Fulman}*{Table 1} for reference.

\begin{table}[H]
    \centering
    \caption{Class Numbers for Exceptional Groups}
    \begin{tabular}{ |c|c|c| }
    \hline
    $G$ & $k(G)\leq$ & Comments \\
    \hline \hline
    $^2B_2(q)$ & $q+3$ & $q=2^{2m+1}$\\ \hline
    $^2G_2(q)$ & $q+8$ & $q=3^{2m+1}$\\ \hline
    $G_2(q)$ & $q^2+2q+9$ &  \\ \hline
    $^2F_4(q)$ & $q^2+4q+17$ & $q=2^{2m+1}$ \\ \hline
    $^3D_4(q)$ & $q^4+q^3+q^2+q+6$ &\\ \hline
    $F_4(q)$ & $q^4+2q^3+7q^2+15q+31$ & \\ \hline
    $E_6(q)$ & $q^6+q^5+2q^4+2q^3+15q^2+21q+60$ &\\ \hline
    $^2E_6(q)$ & $q^6+q^5+2q^4+4q^3+18q^2+26q+62$ & \\ \hline
    $E_7(q)$ & $q^7+q^6+2q^5+7q^4+17q^3+35q^2+71q+103$ & \\ \hline
    $E_8(q)$ & $q^8+q^7+2q^6+3q^5+10q^4+16q^3+40q^2+67q+112$ & \\ \hline
    \end{tabular}
    \label{tbl2}
\end{table}

\begin{enumerate}

\item Let $H \cong E_6(q)$ where $q=p^k$. From the order formula found in \cite{Carter}, $q^{36}$ divides $|E_6(q)|$. From \cite{Bessenrodt}, we know that for any prime $p$, $|n!|_p \leq p^{\frac{n}{p-1}}$. If $q=p^k$, then we have $|n!|_q\leq q^{\frac{n}{k(p-1)}}$ and thus $|\mathrm{A}_n|_q\leq q^{\frac{n}{k(p-1)}}$ where $|\mathrm{A}_n|_p$ is the $p$-part of $\mathrm{A}_n$. By Lemma \ref{order}, $|E_6(q)|$ divides $|\mathrm{A}_n|$ so $q^{36}$ divides $|\mathrm{A}_n|.$ Thus $36 \leq \frac{n}{k(p-1)}$ and $n\geq 36k(p-1)$. Therefore, $|\mathrm{A}_n|\geq|\mathrm{A}_{36k(p-1)}|$.

Now, we note from Table \ref{tbl2} that $|\mathrm{Irr}(E_6(q))|=k(E_6(q))\leq q^6+q^5+2q^4+2q^3+15q^2+21q+60$. Applying Lemma \ref{orderK} gives $|\mathrm{A}_n|<|E_6(q)|\cdot |\mathrm{Irr}(E_6(q))|.$ Hence, $|\mathrm{A}_{36k(p-1)}|<|E_6(q)|\cdot (q^6+q^5+2q^4+2q^3+15q^2+21q+60)$. As with the classical Lie type groups, we can computationally find an upper bound on $p$ and $k$ since the left side grows faster in terms of $p$ and $k$ than the right side. In this case, we find that no values of $p$ and $k$ satisfy the inequality, since substituting $p=2$ and $k=1$ gives $|\mathrm{A}_{36}|>|E_6(2)|\cdot(2^6+2^5+2\cdot 2^4+2\cdot2^3+15\cdot2^2+21\cdot2+60)$. Thus, there are no possible values for $q$ and $n$ such that $\mathrm{cod}(E_6(q))\subseteq \mathrm{cod}(\mathrm{A}_n).$

\item Let $H \cong E_7(q)$ where $q=p^k$. From \cite{Carter}, $q^{63}$ divides $|E_7(q)|$. From Table \ref{tbl2}, $|\mathrm{A}_{63k(p-1)}| < |E_7(q)|\cdot(q^7+q^7+2q^5+7q^4+17q^3+35q^2+71q+103)$. We computationally check that we get a contradiction for $p=2,k=1,$ so there are no possible exceptions.

\item Let $H \cong E_8(q)$ where $q=p^k$. From \cite{Carter}, $q^{120}$ divides $|E_8(q)|$. Thus, using Table \ref{tbl2} as above, we have $|\mathrm{A}_{120k(p-1)}| < |E_8(q)|\cdot(q^8+q^7+2q^6+3q^5+10q^4+16q^3+40q^2+67q+112)$. Now, we computationally check that we get a contradiction for $p=2,k=1,$ so there are no possible exceptions.

\item Let $H \cong F_4(q)$ where $q=p^k$. From \cite{Carter}, $q^{24}$ divides $|F_4(q)|$. From Table \ref{tbl2}, $|\mathrm{A}_{24k(p-1)}| < |F_4(q)|\cdot(q^4+2q^3+7q^2+15q+31)$. Again, we computationally check that we get a contradiction for $p=2,k=1,$ so there are no possible exceptions.

\item Let $H \cong G_2(q)$ where $q=p^k$. From \cite{Carter}, $q^{6}$ divides $|G_2(q)|$. Thus, using Table \ref{tbl2} as above, $|\mathrm{A}_{6k(p-1)}| < |G_2(q)|\cdot(q^2+2q+9)$. Now, we find that $p=2,k=1$ satisfies the inequality, but any other values of $p$ and $k$ do not. However, we note that $G_2(2)$ is not simple, so we instead consider its derived subgroup $G_2(2)'$ (which still satisfies the above inequality). We check for exceptions where $|G_2(2)'|$ divides $|\mathrm{A}_n|$ and $|\mathrm{A}_n|<|G_2(2)'|\cdot(2^2+2\cdot2+9),$ but there are none.

\item Let $H \cong {}^2E_6(q)$ where $q=p^k$. From \cite{Carter}, $q^{36}$ divides $|^2E_6(q)|$. Using Table \ref{tbl2}, $|\mathrm{A}_{36k(p-1)}| < |^2E_6(q)|\cdot(q^6+q^5+2q^4+4q^3+18q^2+26q+62)$. Again, we computationally check that we get a contradiction for $p=2,k=1,$ so there are no possible exceptions.

\item Let $H \cong {}^3D_4(q)$ where $q=p^k$. From \cite{Carter}, $q^{12}$ divides $|^3D_4(q)|$. Thus, using Table \ref{tbl2} similarly to above, $|\mathrm{A}_{12k(p-1)}| < |^3D_4(q)|\cdot(q^4+q^3+q^2+q+6)$. Now, we find that $p=2,k=1$ satisfies the inequality, but any other values of $p$ and $k$ do not. As for the sporadic groups, we check for possible exceptions where $|{}^3D_4(2)|$ divides $|\mathrm{A}_n|$ and $|\mathrm{A}_n|<|{}^3D_4(2)|\cdot(2^4+2^3+2^2+2+2),$ but there are none.

\item Let $H \cong {}^2B_2(q)$ where $q=2^{2m+1}$ and $m\geq 1$. From \cite{Carter}, $q^{2}$ divides $|^2B_2(q)|$. From Table \ref{tbl2}, we have that $|\mathrm{A}_{2(2m+1)}| < |^2B_2(q)|\cdot(q+3)$. In this case, we computationally check that we get a contradiction if $m> 4,$ so $m$ must be less than $5.$ However, checking the divisibility condition, we get no exceptions.

\item Let $H \cong {}^2F_4(q)$ where $q=2^{2m+1}$ and $m\geq 1$. From \cite{Carter}, $q^{12}$ divides $|^2F_4(q)|$. Thus, using Table \ref{tbl2} as above, $|\mathrm{A}_{12(2m+1)}| < |^2F_4(q)|\cdot(q^2+4q+17)$. Now, we computationally check that we get a contradiction for $m = 1,$ so there are no exceptions

\item Let $H \cong {}^2G_2(q)$ where $q=3^{2m+1}$ and $m\geq 1$. From \cite{Carter}, $q^{3}$ divides $|^2G_2(q)|$. From Table \ref{tbl2}, $|\mathrm{A}_{3(2m+1)\cdot2}| < |^2G_2(q)|\cdot(q+8)$. Again, we computationally check that we get a contradiction for $m= 1,$ so there are no exceptions.
\end{enumerate}
\end{proof}

\begin{Theorem}\label{GmodN}
    Let $G$ be a finite group such that $\cod(G)=\cod(\mathrm{A}_n)$ where $n\geq 5$. Let $N$ be a maximal subgroup of $G$. Then, $G/N\cong \mathrm{A}_n$.
\end{Theorem}
\begin{proof}
    By Lemma \ref{perfect}, $G$ is perfect. Thus $G/N$ is a nonabelian simple group. By Lemma \ref{subset}, we have $\mathrm{cod}(G/N) \subseteq \mathrm{cod}(G)=\mathrm{cod}(\mathrm{A}_n)$. By Lemmas \ref{alternating}, \ref{sporadic}, \ref{classical}, and \ref{exceptional}, $G/N$ cannot be an alternating group of degree $m\neq n$, a sporadic simple group or the Tits group, a classical simple group of Lie type, or an exceptional simple group of Lie type. Thus, $G/N \cong \mathrm{A}_n$.
\end{proof}

Now we present the proof of Theorem \ref{thm11}.
\begin{proof}

Let $G$ be a minimal counterexample and $N$ be a maximal normal subgroup of $G$. By Lemma \ref{perfect}, $G$ is perfect, and by Theorem \ref{GmodN}, $G/N \cong \mathrm{A}_n$. In particular, $N \neq 1$ as $G \not\cong \mathrm{A}_n$.

{\bf Step 1:} $N$ is a minimal normal subgroup of $G$.

Suppose $L$ is a non-trivial normal subgroup of $G$ with $L < N$. Then by Lemma \ref{subset}, we have $\mathrm{cod}(G/N) \subseteq \mathrm{cod}(G/L) \subseteq \mathrm{cod}(G)$. However, $\mathrm{cod}(G/N)=\mathrm{cod}(\mathrm{A}_n)=\mathrm{cod}(G)$ so equality must be obtained in each inclusion. Thus, $\mathrm{cod}(G/L)=\mathrm{cod}(\mathrm{A}_n)$ which implies that $G/L \cong \mathrm{A}_n$ since $G$ is a minimal counterexample. This is a contradiction since we also have $G/N\cong \mathrm{A}_n,$ but $L < N$.

{\bf Step 2:} $N$ is the only non-trivial, proper normal subgroup of $G$.

Otherwise we assume $M$ is another proper nontrivial  normal subgroup of $G$. If $N$ is included in $M$, then $M=N$ or $M=G$ since $G/N$ is simple, a contradiction. Then $N\cap M=1$ and $G=N\times M$. Since  $M$ is also a maximal normal subgroup of $G$, we have $N\cong M\cong \mathrm{A}_n$. Choose $\psi_1\in \irr(N)$ and $\psi_2\in \irr(M)$ such that $\cod(\psi_1)=\cod(\psi_2)=\max(\cod(\mathrm{A}_n))$. Set $\chi=\psi_1\cdot\psi_2\in \irr(G)$. Then $\cod(\chi)=(\max(\cod(\mathrm{A}_n)))^2\notin \cod(G)$, a contradiction.


{\bf Step 3:} For each non-trivial $\chi \in \mathrm{Irr}(G|N):=\mathrm{Irr}(G)-\mathrm{Irr}(G/N), \chi$ is faithful.

We construct $\mathrm{Irr}(G/N)$ as the same as Lemma \ref{codegreeSubset}. Then it follows by the definition of $\mathrm{Irr}(G|N)$ that if $\chi \in \mathrm{Irr}(G|N),$ $N \not\leq \mathrm{ker}(\chi).$ Thus since $N$ is the unique nontrivial, proper, normal subgroup of $G$, $\mathrm{ker}(\chi) = G$ or $\mathrm{ker}(\chi) = 1$. Therefore, $\mathrm{ker}(\chi) = 1$ for all nontrivial $\chi \in \mathrm{Irr}(G|N).$

{\bf Step 4:} $N$ is an elementary abelian group.

Suppose that $N$ is not abelian. Since $N$ is a minimal normal subgroup, by \cite{Dixon}*{Theorem 4.3A (iii)}, $N=S^n$ where $S$ is a nonabelian simple group and $n\in \mathbb{Z}^+$. By Lemmas \ref{oneextend} and \ref{prodsimple}, there is a non-trivial character $\chi\in \mathrm{Irr}(N)$ which extends to some $\psi\in\mathrm{Irr}(G).$ Now, ker$(\psi)=1$ by Step 3, so cod$(\psi)=|G|/\psi(1)=|G/N|\cdot |N|/\chi(1).$ However, by assumption, we have that $\mathrm{cod}(G)= \mathrm{cod}(\mathrm{A}_n)= \mathrm{cod}(G/N)$. Thus, cod$(\psi) \in \mathrm{cod}(G)= \mathrm{cod}(G/N),$ so $\mathrm{cod}(\psi)=|G/N|/\phi(1)$ for some $\phi\in\textrm{Irr}(G/N).$ Hence, $|G/N|$ is divisible by $\mathrm{cod}(\psi)$ which contradicts the fact that $\mathrm{cod}(\psi)=|G/N|\cdot |N|/\chi(1),$ as $\chi(1)\neq |N|.$ Thus $N$ must be abelian.

Now to show that $N$ is elementary abelian, let a prime $p$ divide $|N|.$ Then $N$ has a $p$-Sylow subgroup $K$, and $K$ is the unique $p$-Sylow subgroup of $N$ since $N$ is abelian, so $K$ is characteristic in $N$. Thus, $K$ is a normal subgroup of $G,$ so $K=N$ as $N$ is minimal. Thus $|N|=p^n.$ Now, take the subgroup $N^p=\{n^p \mid n \in N\}$ of $N,$ which is proper by Cauchy's theorem. Since $N^p$ is characteristic in $N,$ it must be normal in $G,$ so $N^p$ is trivial by the uniqueness of $N.$ Thus every element of $N$ has order $p,$ and $N$ is elementary abelian.

{\bf Step 5:} $\mathbf{C}_G(N) = N.$

First note that since $N$ is normal, $\mathbf{C}_G(N) \trianglelefteq G.$ Additionally, since $N$ is abelian by Step 4, $N \leq \mathbf{C}_G(N)$. By the maximality of $N,$ we must have $\mathbf{C}_G(N) = N$ or $\mathbf{C}_G(N) = G.$ If $\mathbf{C}_G(N) = N,$ we are done.

If not, then $\mathbf{C}_G(N) = G,$ so $N$ must be in the center of $G.$ Then since $N$ is the unique minimal normal subgroup of $G$ by Step 2, we must have that $|N|$ is prime. If not, there always exists a proper non-trivial subgroup $K$ of $N,$ and $K$ is normal since it is contained in $\mathbf{Z}(G),$ contradicting the minimality of $N.$ Moreover, since $G$ is perfect, we have that $\mathbf{Z}(G) = N,$ and $N$ is isomorphic to a subgroup of the Schur multiplier of $G/N$ \cite{Isaacs}*{Corollary 11.20}.

Now, we note that it is well-known that for $n>7$, the Schur multiplier of $\mathrm{A}_n$ is $\mathbb{Z}_2,$ so $G\cong 2.\mathrm{A}_n$. From \cite{Malle}, $2.\mathrm{A}_n$ always has a character degree of order $2^{\lfloor (n-2)/2\rfloor}$. Let $\chi$ be such an irreducible character of $2.\mathrm{A}_n$ with $\chi(1)=2^{\lfloor (n-2)/2\rfloor}.$ Recall that by Step 2, there is only one non-trivial proper normal subgroup of $G \cong 2.\mathrm{A}_n$. In particular $N \cong \mathbb{Z}_2$ is the only non-trivial proper normal subgroup of $G$. Thus $|\mathrm{ker}(\chi)|=1$ or $2$. Then we have $\cod(\chi)=\frac{|2.\mathrm{A}_n:\mathrm{ker}(\chi)|}{\chi(1)}$. If $|\mathrm{ker}(\chi)|=1,$ then $\cod(\chi)=\frac{n!}{2^{\lfloor (n-2)/2\rfloor}},$ and if $|\mathrm{ker}(\chi)|=2,$ then $\cod(\chi)=\frac{n!/2}{2^{\lfloor (n-2)/2\rfloor}}= \frac{n!}{2^{\lfloor n/2\rfloor}}.$ In either case, for any prime $p\neq 2,|\cod(\chi)|_p=|n!|_p=|\mathrm{A}_n|_p$. Since $\cod(G)=\cod(\mathrm{A}_n),$ we know that $\cod(\chi)\in \cod(\mathrm{A}_n)$. Therefore, there is a character degree of $\mathrm{A}_n$ which is a power of $2$.

However, from \cite{Malle}, we know that for $n>7, \mathrm{A}_n$ only has a character degree equal to a power of $2$ when $n=2^d+1$ for some positive integer $d$. In this case, $2^d=n-1\in \mathrm{cd}(\mathrm{A}_n)$ so we need $\frac{|\mathrm{A}_n|}{n-1}=\frac{|2.\mathrm{A}_n|}{2^{\lfloor (n-2)/2\rfloor}}$ or $\frac{|2.\mathrm{A}_n|}{2^{\lfloor n/2\rfloor}}$. Hence, $\frac{1}{n-1}=\frac{2}{2^{\lfloor (n-2)/2\rfloor}}=\frac{1}{2^{\lfloor (n-2)/2\rfloor-1}}$ or $\frac{1}{2^{\lfloor n/2\rfloor-1}}$ so $n-1=2^{\lfloor (n-2)/2\rfloor-1}$ or $2^{\lfloor n/2\rfloor-1}$. However, the only integer solution to either of these equations occurs when $n=9$ and  $9-1=8=2^3=2^{\lfloor 9/2\rfloor-1}$. In this case, we check the ATLAS \cite{Atlas} to find that the codegree sets of $A_9$ and $2.A_9$ do not have the same order. This is a contradiction, so $\mathbf{C}_G(N)=N$.

{\bf Step 6:} Let $\lambda$ be a non-trivial character in $\mathrm{Irr}(N)$ and $\vartheta \in \mathrm{Irr}(I_G(\lambda)|\lambda),$ the set of irreducible constituents of $\lambda^{I_G(\lambda)},$ where $I_G(\lambda)$ is the inertia group of $\lambda \in G.$ Then $\frac{|I_G(\lambda)|}{\vartheta(1)} \in \mathrm{cod}(G).$ Also, $\vartheta(1)$ divides $|I_G(\lambda)/N|,$ and $|N|$ divides $|G/N|.$ Lastly, $I_G(\lambda) < G,$ i.e. $\lambda$ is not $G$-invariant.

Let $\lambda$ be a non-trivial character in $\operatorname{Irr}(N)$ and $\vartheta \in \operatorname{Irr}(I_G(\lambda)|\lambda)$. Let $\chi$ be an irreducible constituent of $\vartheta^G.$ By \cite{Isaacs}*{Corollary 5.4}, we know $\chi \in \operatorname{Irr}(G)$, and by \cite{Isaacs}*{Definition 5.1}, we have $\chi(1) = \frac{|G|}{|I_G(\lambda)|} \cdot \vartheta(1)$. Moreover, we know tat $\operatorname{ker}(\chi) = 1$ by Step 2, and thus $\cod(\chi) = \frac{|G|}{\chi(1)} = \frac{|I_G(\lambda)|}{\vartheta(1)}$, so $\frac{|I_G(\lambda)|}{\vartheta(1)} \in \cod(G)$. Now, since $N$ is abelian, $\lambda(1) = 1$, so we have $\vartheta(1) = \vartheta(1)/\lambda(1)$ which divides $\frac{|I_G(\lambda)|}{|N|}$, so $|N|$ divides $\frac{|I_G(\lambda)|}{\vartheta(1)}$. Moreover, we know that $\cod(G) = \cod(G/N),$ and all elements in $\cod(G/N)$ divide $|G/N|$, so $|N|$ divides $|G/N|$.

Next, we want to show $I_G(\lambda)$ is a proper subgroup of $G$. To reach a contradiction, assume $I_G(\lambda) = G$. Then $\operatorname{ker}(\lambda) \unlhd G$. From Step 2, we know $\operatorname{ker}(\lambda) = 1,$ and from Step 4, we know $N$ is a cyclic group of prime order. Thus by the Normalizer-Centralizer theorem, we have $G / N= \mathbf{N}_{G}(N) / \mathbf{C}_{G}(N) \leq \operatorname{Aut}(N)$ so $G / N$ is abelian, a contradiction.

{\bf Step 7:} Final contradiction.

From Step 4, $N$ is an elementary abelian group of order $p^m$ for some prime $p$ and integer $m\geq1$. By the Normalizer-Centralizer theorem, $\mathrm{A}_n \cong G/N = \mathbf{N}_G(N)/\mathbf{C}_G(N) \leq \mathrm{Aut}(N)$ and $m>1$. Note that in general, $\mathrm{Aut}(N)=\mathrm{GL}(m,p)$. By Step 6, $|N|$ divides $|G/N|,$ so we know that $|N|=p^m$ divides $|\mathrm{A}_n|$ and $G/N\cong \mathrm{A}_n\lesssim \mathrm{GL}(m,p).$ We prove by contradiction that this cannot occur.

First, we claim that if $p^m$ divides $|\mathrm{A}_n|$ and $\mathrm{A}_n\lesssim (\mathrm{GL}(m,p),$ then $p$ must equal $2$. To show this, we note that for $p>2,$ by \cite{Bessenrodt}, we have that if $p^m$ divides $|\mathrm{A}_n|,$ then $m < \frac{n}{2}$. However, Theorem 1.1 of \cite{Wagner1} shows that if $n>6$, the minimal faithful degree of a modular representation of $\mathrm{A}_n$ over a field of characteristic $p$ is at least $n-2$. Since embedding $\mathrm{A}_n$ as a subgroup of $\mathrm{GL}(m,p)$ is equivalent to giving a faithful representation of degree $m$ over a field of characteristic $p$, we have that $m\geq n-2$. This is a contradiction since $\frac{n}{2}> n-2$ implies $n<4$. Therefore, $p=2$.

Now, let $p=2$. As above, from \cite{Bessenrodt}, we obtain $|n!|_2\leq 2^{n-1}$. Thus, if $2^m$ divides $|\mathrm{A}_n|,$ then $m \leq |\mathrm{A}_n|_2 \leq 2^{n-2}$. Now, Theorem 1.1 of \cite{Wagner2} shows that if $n>8$, then the minimal faithful degree of a modular representation of $\mathrm{A}_n$ over a field of characteristic $2$ is at least $n-2$. Therefore, we must have $m\geq n-2$, so $m=|\mathrm{A}_n|_2=2^{n-2}$ is the only option.

Let $\lambda \in \mathrm{
Irr}(N),$ $\vartheta \in \mathrm{Irr}(I_G(\lambda)|\lambda),$ and $T := I_G(\lambda)$. Then $1<|G:T|<|N|=2^{n-2}$ for $|G:T|$ is the number of all conjugates of $\lambda$. By Step 5, we know that $\frac{|T|}{\vartheta(1)}\in\mathrm{cod}(G)$ and moreover that $|N|$ divides $\frac{|T|}{\vartheta(1)}$. Since $|N|_2=|\mathrm{A}_n|_2$ and $\mathrm{cod}(G)=\mathrm{cod}(\mathrm{A}_n),$ we know that $\left|\frac{|T|}{\vartheta(1)}\right|_2=|N|_2.$ Thus $\left|\frac{|T/N|}{\vartheta(1)}\right|_2=1$ so the $2$-parts of $|T/N|$ and $\vartheta(1)$ are equal. Thus for every $\vartheta \in \text{Irr}(T\mid \lambda)$, we have $|\vartheta(1)|_2=|T/N|_2.$ However, $|T/N|= \sum_{\vartheta\in\text{Irr}(T\mid\lambda)} \vartheta(1)^2.$ Hence, if $|\vartheta(1)|_2=2^k\geq 2$ for every $\vartheta\in \text{Irr}(T\mid\lambda),$ we would have $|T/N|_2=2^{2k}$ contradicting the fact that $|\vartheta(1)|_2=|T/N|_2.$ Therefore, $|T/N|_2=1.$ Thus, since $|G/N|_2\geq|N|_2=2^{n-2}$, we have $|G:T|_2=|G/N:T/N|_2\geq 2^{n-2},$ so $|G:T|\geq 2^{n-2}=|N|,$ which is a contradiction.

We have one final exception to consider: $n=8$, $p=2,$ and $m=4,5$ or $6$. In this case, $\mathrm{A}_8 \cong \mathrm{GL}(4,2)$ and $2^6$ divides $|\mathrm{A}_8|$. Now, $\mathrm{cod}(\mathrm{A}_8)=\{1, 2^6\cdot3^2 \cdot5, 2^5 \cdot3^2 \cdot5, 2^4 \cdot3^2 \cdot7, 2^6 \cdot3\cdot5, 2^4 \cdot3^2 \cdot5, 2^6 \cdot3^2, 2^6 \cdot7, 2^3 \cdot3^2 \cdot5, 3^2 \cdot5\cdot7, 2^5 \cdot3^2\}$ from \cite{Gintz}. We will look at each possibility for $m$ in turn.

First, let $m=4$. Then we have $G/N \cong \text{A}_8 \cong \text{GL}(4,2), N=(\mathbb{Z}_2)^4$ so $G$ is an extension of $\text{GL}(4,2)$ by $N.$ Suppose first that this extension is split and $G$ is a semidirect product. This semidirect product is defined by a homomorphism $\phi: \text{GL}(4,2)\to \text{Aut}((\mathbb{Z}_2)^4)\cong \text{GL}(4,2).$ However, since $\text{GL}(4,2)$ is simple, $\text{ker}(\phi)=1$ or $\text{GL}(4,2).$ In the first case, we have the trivial direct product, so there are at least two copies of $\text{GL}(4,2)$ as normal subgroups of $G,$ which contradicts Step 2. In the second case, $\phi$ is some automorphism of $\text{GL}(4,2).$ Here, we can check using GAP that any such $\phi$ creates a semidirect product $\mathrm{GL}(4,2)\rtimes_\phi (\mathbb{Z}_2)^4$ which does not have the same codegree set as $\mathrm{A}_8.$ Now, suppose that the extension is non-split. Then, \cite{Basheer} gives that there is a unique non-split extension $2^4.\text{GL}(4,2).$ However, we find using GAP that it doesn't have the same codegree set as $\mathrm{A}_8.$


Second, let $m=5$. As above, $|G:T|<|N|=2^5$ and $\frac{|T|}{\vartheta(1)}\in\mathrm{cod}(G)$ such that $2^5$ divides $\frac{|T|}{\vartheta(1)}$. Further, $|\frac{|T/N|}{\vartheta(1)}|_2\leq2$ so $|T/N|_2\leq 4$ and $|G/N:T/N|_2\geq 16$. Thus, we have $16$ divides $|G/N:T/N|$ and $|G/N:T/N|< 32.$ But we check the index of all subgroups of $G/N\cong\text{A}_8$ using GAP and find that none of them satisfy these two properties.

Finally, let $m=6$. Now, $|N|_2=|\mathrm{A}_8|_2.$ For this case the same argument as above for general $\mathrm{A}_n$ holds, and we reach a contradiction. Thus we find that every $|N|=p^m$ produces a contradiction, so $N=1$ and $G\cong \mathrm{A}_n$.
\end{proof}

\end{section}

\begin{section}{Acknowledgements}

This research was conducted under NSF-REU grant DMS-1757233, DMS-2150205 and NSA grant H98230-21-1-0333, H98230-22-1-0022 by Dolorfino, Martin, Slonim, and Sun during the Summer of 2022 under the supervision of Yang. The authors gratefully acknowledge the financial support of NSF and NSA, and also thank Texas State University for providing a great working environment and support. Yang was also partially supported by grants from the Simons Foundation (\#499532, \#918096, to YY). The authors would also like to thank Prof. Richard Stanley for his help.
\end{section}


\begin{thebibliography}{99}


    \bibitem{Ahanjideh} N. Ahanjideh, Nondivisibility among irreducible character co-degrees. \textit{Bull. Aust. Math. Soc.}, \textbf{105} (2022), 68-74.


    \bibitem{Aziziheris} K. Aziziheris, F. Shafiei, F. Shirjian, Simple groups with few irreducible character degrees. \textit{J. Algebra Appl.}, \textbf{20} (2021), 2150139.

    \bibitem{Bahri} A. Bahri, Z. Akhlaghi, B. Khosravi, An analogue of Huppert's conjecture for character codegrees. \textit{Bull. Aust. Math. Soc.}, \textbf{104} (2021), no. 2, 278-286.

    \bibitem{Basheer} A. B. M. Basheer and J. Moori, Fischer Matrices of Dempwolff Group $2^5.\text{GL}(5,2).$ \textit{Int. J. Group Theory}, \textbf{1} (2012), 43-63.

    \bibitem{Bessenrodt} C. Bessenrodt, H. P. Tong-Viet, J. Zhang, Huppert's conjecture for alternating groups. \textit{J. Algebra}, \textbf{470} (2017), 353-378.

    \bibitem{Julia} J. Bezanson, S. Karpinski, V. B. Shah, A. Edelman, Julia: A fast dynamic language for technical computing. \textit{ArXiv Preprint}, ArXiv:1209.5145.

    \bibitem{Carter} R. W. Carter, \textit{Simple Groups of Lie Type.} Wiley, 1989.

    \bibitem{Chillag} D. Chillag and M. Herzog, On character degrees quotients. \textit{Arch. Math.}, \textbf{55} (1990), 25-29.

    \bibitem{Atlas} J. H. Conway et. al, \textit{Atlas of Finite Groups.} Oxford Clarendon Press, 1985.

    \bibitem{Dixon} J. D. Dixon and B. Mortimer, \textit{Permutation Groups.} Spring, 1996.

    \bibitem{Dolorfino} M. Dolorfino, L. Martin, Z. Slonim, Y. Sun, Y. Yang, On the characterization of sporadic simple groups by codegrees. submitted.

    \bibitem{Fulman} J. Fulman and R. Guralnick, Bounds on the number and sizes of conjugacy classes in finite Chevalley groups with applications to derangements. \textit{Trans. Amer. Math. Soc.}, \textbf{364} (2012), 3023-3070.


    \bibitem{Gintz} M. Gintz, M. Kortje, M. laurence, Y. Liu, Z. Wang, Y. Yang, On the characterization of some nonabelian simple groups with few codegrees. \textit{Comm. Algebra}, \textbf{50} (2022), 3932-3939.

    \bibitem{gzy} H. Guan, X. Zhang, Y. Yang, Recognizing Ree groups ${}^2G_2(q)$ using the codegree set. \textit{Bull. Aust. Math. Soc.}, \url{https://www.doi.org/10.1017/S0004972722001022}.



    \bibitem{Hung} N. N. Hung, Group pseudo-algebras of finite simple groups. {In progress}.

    \bibitem{Isaacs} I. M. Isaacs, \textit{Character Theory of Finite Groups.} New York Academic Press, 1976.

    \bibitem{James} G. James and A. Kerber, \textit{The Representation Theory of the Symmetric Group.} Addison-Wesley Publishing Company, 1981.


 \bibitem{Khukhro} E. I. Khukrho and V. D. Mazurov, \textit{Unsolved Problems in Group Theory. The Kourovka Notebook. No. 20.} Russian Academy of Sciences, 2022.

    \bibitem{lmy}  Y. Liu and Y. Yang, Huppert's analogue conjecture for $\PSL(3,q)$ and $\PSU(3,q)$. \textit{Results Math.}, \textbf{78} (2023), No. 7.

    \bibitem{Malle} G. Malle and A.E. Zalesskii, Prime power degree representations of quasi-simple groups. \textit{Arch. Math.}, \textbf{77} (2001), 461-468.

    \bibitem{Moreto} A. Moret\'o, Complex group algebra of finite groups: Brauer's problem 1. \textit{Adv. Math.}, \textbf{208} (2007), 236-248.

    \bibitem{Qian} G. Qian, Y. Wang, H. Wei, Co-degrees of irreducible characters in finite groups. \textit{J. Algebra}, \textbf{312} (2007), 946-955.

    \bibitem{Wagner2} A. Wagner, The faithful linear representations of least degree of $\rm{S_n}$ and $\rm{A_n}$ over a field of characteristic 2. \textit{Math. Z.}, \textbf{151} (1976), 127-138.

    \bibitem{Wagner1} A. Wagner, The faithful linear representations of least degree of $\rm{S_n}$ and $\rm{A_n}$ over a field of odd characteristics. \textit{Math. Z.}, \textbf{154} (1977), 104-113.

\end{thebibliography}
\end{document}